\documentclass[a4paper,12pt]{article}
\usepackage{graphicx}
\usepackage{cite}        
\usepackage{amssymb,amsmath,latexsym,enumerate}
\usepackage{subfigure,epsfig}
\usepackage{array,longtable,lscape,bigints}
\usepackage{color}
\usepackage{enumerate}
\usepackage{authblk}

\begin{document}

\title{Lax representation and quadratic first integrals for a family of non-autonomous second-order differential equations}

\author[1,2]{Dmitry I. Sinelshchikov}
\author[2]{Ilia Yu. Gaiur}
\author[2]{Nikolay A. Kudryashov}
\affil[1]{Department of Applied Mathematics, National Research University Higher School of Economics, 34 Tallinskaya street, 123458, Moscow, Russian Federation}
\affil[2]{Department of Applied Mathematics, National Research Nuclear University MEPhI,  31 Kashirskoe Shosse, 115409 Moscow, Russian Federation}

\maketitle

\begin{abstract}
We consider a family of non-autonomous second-order differential equations, which generalizes the Li\'enard equation. We explicitly find the necessary and sufficient conditions for members of this family of equations to admit quadratic, with the respect to the first derivative, first integrals. We show that these conditions are equivalent to the conditions for equations in the family under consideration to possess Lax representations. This provides a connection between the existence of a quadratic first integral and a Lax representation for the studied dissipative differential equations, which may be considered as an analogue to the theorem that connects Lax integrability and Arnold--Liouville integrability of Hamiltonian systems. We illustrate our results by several examples of dissipative equations, including generalizations of the Van der Pol and Duffing equations, each of which have both a quadratic first integral and a Lax representation.
\end{abstract}

\section{Introduction}
We consider the following family of equations
\begin{equation}\label{eq0.1}
  y_{zz}+f(z, y)y_z+g(z, y)=0,
\end{equation}
which generalizes the Li\'enard equation \cite{Polyanin}. Here $f(z, y)$ and $g(z, y)$ are arbitrary functions, which do not vanish, i.e., $f\neq0$ and $g\neq0$. We also exclude from the consideration the case of a linear equation by assuming that $|f_{y}|^{2}+|g_{yy}|^{2}\neq0$.

In the case of $f_{z}=g_{z}=0$ from (\ref{eq0.1}) we obtain the classical Li\'enard equation
\begin{equation}\label{eq0.2}
   y_{zz}+f(y)y_z+g(y)=0.
\end{equation}
Further we study both the general case of \eqref{eq0.1} and the classical Li\'enard equation (\ref{eq0.2}).

Equations of type (\ref{eq0.1}) play an important role in various fields of science, such as nonlinear dynamics, mechanics and biology (see, e.g., \cite{Polyanin,Guckenheimer1983}). In addition, equations from family (\ref{eq0.1}) often appear as symmetry reductions of second- and third-order nonlinear partial differential equations.

Integrability of family (\ref{eq0.1}) and some of its particular cases has been considered in various works. Perhaps the most well-studied equation from (\ref{eq0.1}) is the classical Li\'enard equation. For example, Lie and Noether symmetries of (\ref{eq0.2}) were studied in \cite{Pandey2009,Pandey2009a,Gubbiotti2014,Nucci2011}. Linearization and equivalence to some Painlev\'e--Gambier equations via the generalized Sundman transformations were considered in \cite{Duarte1994,Nakpim2010a,Moyo2011,Euler2004,Kudryashov2016, Kudryashov2016a, Kudryashov2017a}. In works \cite{DAmbrosi2009,Nucci2010,Nucci2010a,Kudryashov2017,Sinelshchikov2017} Lagrangians and Jacobi last multipliers for the Li\'enard-type equations were studied. As far as integrability of (\ref{eq0.1}) is concerned, its linearizability conditions via both point and Sundman transformations can be obtained from the results of \cite{Lie1888,Tresse1896,MAHOMED1989}. Notice also that complete classification of Lie symmetries of second-order differential equations, at most cubic in the first derivative, was obtained in \cite{Bagderina2013}. In \cite{MURIEL2009,Muriel2011,Meleshko2013,Bagderina2016} linear and certain rational first integrals for (\ref{eq0.1}) and its generalization were classified. Connections of equations from family (\ref{eq0.1}) with some Painlev\'e--Gambier equations via nonlocal transformations were considered in \cite{Sinelshchikov2018}. However, to the best of our knowledge, neither Lax integrability nor classification of quadratic first integrals for (\ref{eq0.1}) have been studied previously.

Therefore, the main aim of this work is to study equation (\ref{eq0.1}) as a Lax integrable system and obtain the necessary and sufficient conditions for equation (\ref{eq0.1}) to possess a first integral, which is quadratic with respect to the first derivative. Lax integrability requires a differential system to have a Lax representation, which also called the Lax pair. There are many examples of Hamiltonian systems that are Lax-integrable systems, for instance, the Calogero--Moser system, the Toda chain and so on~\cite{Perelomov1990,Babelon2003}. In this work we generalize the Lax pair for the harmonic oscillator to find equations of type (\ref{eq0.1}), which has a Lax representation. The knowledge of a Lax pair for an equation from family (\ref{eq0.1}) leads to the corresponding quadratic first integral. On the other hand, one can consider the problem of finding all equations from family (\ref{eq0.1}), which admit quadratic first integrals. Below, we solve this problem and obtain explicit correlations on functions $f$ and $g$, which provide the necessary and sufficient conditions for the existence of a quadratic first integral. Finally, we show that these conditions are equivalent to conditions for the existence of a Lax representation for an equation from family (\ref{eq0.1}).

The rest of this work is organized as follows. In Section 2 we show that equation (\ref{eq0.1}) admits a Lax representation provided certain correlations on functions $f$ and $g$ hold. Section 3 is devoted to establishing necessary and sufficient conditions for the existence of a quadratic first integral for (\ref{eq0.1}). We provide several examples of Lax integrable equations from family (\ref{eq0.1}) in Section 4. In the last section we show that the conditions for the existence of a Lax representation and a quadratic first integral are equivalent and briefly discuss our results.

\section{Lax representation}

In this section we give basic notion of the Lax integrability and present Lax pairs corresponding to (\ref{eq0.1}) and \eqref{eq0.2}. Let us note that the statements presented below concerning Lax integrability are true not only for Hamiltonian systems, but also for any system of differential equations.

By definition, a Lax integrable system can be presented as an isospectral deformation of the following linear problem \cite{Babelon2003, Torrielli2016}
\begin{equation}\label{eq0.4}
\begin{gathered}
  \left(L - \lambda\right) \Psi = 0, \\
  \left(\partial_z + M \right)\Psi = 0,
\end{gathered}
\end{equation}
where $L$ and $M$ are $2\times2$ complex matrices which entries depend on $z$ and $\Psi$ is a complex two-dimensional vector function of $z$. The isospectrality condition ($\lambda_z=0$) gives the Lax equation
\begin{equation}\label{eq0.2.2}
 L_z = [L,M],
\end{equation}
which is the compatibility condition for (\ref{eq0.4}). Lax equation (\ref{eq0.2.2}) provides nontrivial conditions on elements of the $L$ matrix, which give us the corresponding integrable system.
It is easy to show, that if equation (\ref{eq0.2.2}) holds, then \cite{Babelon2003, Torrielli2016}
\begin{equation}\label{eq0.3.1}
  (L^k)_z = [L^k, M].
\end{equation}

Since $\operatorname{tr}(AB)=\operatorname{tr}(BA)$, where $\operatorname{tr}$ denotes the trace of a square matrix, and
\begin{displaymath}
\operatorname{tr}(L_z)=\frac{{\rm d}}{{\rm d}z}\operatorname{tr}(L),
\end{displaymath}
the Lax equation allows us to find first integrals of the corresponding integrable system in the following form
\begin{equation}\label{eq0.4.1}
  I_k = \frac{\operatorname{tr}(L^k)}{k!}, \quad \frac{{\rm d}}{{\rm d}z}{I_k} = 0.
\end{equation}
On the other hand, since the Lax equation is obtained from isospectral deformation ($\lambda_z = 0)$, first integrals can be found as eigenvalues of $L$.

Now we generalize the $L$ matrix for the harmonic oscillator \cite{Babelon2003, Torrielli2016} in order to obtain equation (\ref{eq0.1}) as an isospectral condition for linear system (\ref{eq0.4}). Namely, we consider the $L$ matrix as an element of the algebra $\mathfrak{sl}(2, \mathbb{C})$, i.e.
\begin{equation}\label{0.3.1}
  L = \left(\begin{array}{cc}
              y_z + F & U \\
              U & -(y_z + F)
            \end{array}\right),
\end{equation}
where $F=F(z,y)$ and $U=U(z,y)$ are arbitrary nonzero functions. We do not impose any constraints on the $M$ matrix until we write it explicitly for some particular cases of equation (\ref{eq0.1}). With the help of $L$ matrix (\ref{0.3.1}), we obtain the corresponding first integral for equation (\ref{eq0.1})
\begin{equation}\label{0.5}
  \mathcal{I}= \left(y_z + F\right)^2 + U^2.
\end{equation}

Let us present conditions on $f$ and $g$ for equation (\ref{eq0.1}) to possess a Lax representation and the corresponding Lax pairs for (\ref{eq0.1}).

\subsection{Generic case of \eqref{eq0.1}}

Substituting (\ref{0.3.1}) into (\ref{eq0.2.2}) and taking into account (\ref{eq0.1}) we obtain that if coefficients $f$ and $g$ satisfy the conditions
\begin{equation}\label{eq3.1.1.2.1}
2g = \left(U^2\right)_{y} + 2F_z ,\quad f = F_y,\quad U_z-F U_y=0,
\end{equation}
then equation (\ref{eq0.1}) admits the Lax pair
\begin{equation}\label{eq3.3}
  L=\left(\begin{array}{cc}
            y_z + F  & U \\
           U & - y_z - F
          \end{array}\right), \quad
          M = \left(\begin{array}{cc}
                0 & \frac{1}{2}U_y \\
                -\frac{1}{2}U_y & 0
              \end{array}\right).
\end{equation}
Correlations \eqref{eq3.1.1.2.1} provide an overdetermined system of equations for the functions $F$ and $U$. The corresponding compatibility conditions give us the necessary and sufficient conditions for the existence of a Lax representation with the $L$ matrix of form \eqref{0.3.1} for equation (\ref{eq0.1}). Let us remark, that these compatibility conditions will be found in the explicit form in the next section.

\subsection{The classical Li\'enard equation}
Let us consider the classical Li\'enard equation, i.e. we assume that $g_z=0$ and $f_z=0$. In this case, from (\ref{eq3.1.1.2.1}) one can find the expressions for $U$ and $F$ as follows
\begin{equation}\label{2.2}
\begin{gathered}
   U^2 = 2 \int\limits^y_{y_0}{\rm d}q\left(\kappa\int\limits^q_{q_0} f(\xi){\rm d}\xi + \mu\right)^{-1} + 2 \kappa^{-1}z, \\
   F = \int\limits^y_{y_0} f(\xi){\rm d}\xi + \mu,
\end{gathered}
\end{equation}
where $\mu$ and $\kappa\neq0$ are arbitrary constants and here and below $\xi$ and $q$ are dummy integration variables.

Since $2g=\left(U^2\right)_y$, we obtain that the family of classical Li\'enard equations
\begin{equation}\label{2.3}
  y_{zz} + f y_z + \left(\kappa\int\limits^y_{y_0} f(\xi){\rm d}\xi + \mu\right)^{-1} = 0,
\end{equation}
has the Lax pair
\begin{equation}\label{2.2.1}
  L=\left(\begin{array}{cc}
            y_t+\int\limits^y_{y_0} f(\xi){\rm d}\xi + \mu & U \\
            U  & -y_t-\int\limits^y_{y_0} f(\xi){\rm d}\xi - \mu
          \end{array}\right),
\end{equation}
where $U$ is defined by (\ref{2.2}) and $M$-matrix can be found from (\ref{eq3.3}). Consequently, family of equations (\ref{2.3}) admits the first integral
\begin{equation}\label{2.4}
  \mathcal{I} = \left(y_z+\int\limits^y_{y_0} f(\xi){\rm d}\xi + \mu \right)^2 + 2 \int\limits^y_{y_0}{\rm d}q\left(\kappa\int\limits^q_{q_0} f(\xi){\rm d}\xi + \mu\right)^{-1} + 2 \kappa^{-1}z.
\end{equation}

On the other hand, one can find the value of $f$ for a given $g$, which is
\begin{equation}\label{2.5}
  f= - \nu g^{-2} g_{y},
\end{equation}
where $\nu\neq 0$ is an arbitrary constant. Rewriting \eqref{2.4} in terms of $g$ we obtain
\begin{equation}
   \mathcal{I} = \left(y_z + \frac{\nu}{g}\right)^2 + 2\left( \int\limits^y_{y_0}g(\xi){\rm d}\xi + \nu z\right).
\end{equation}
Therefore, we have explicitly found a family of the classical Li\'enard equations, defined either by \eqref{2.3} or \eqref{2.5}, with an energy-like first integral and a Lax representation.

\section{Necessary and sufficient conditions for the existence of a quadratic first integral}

We have shown above that equations of type (\ref{eq0.1}) can possess a Lax representation, provided certain correlations on $f$ and $g$ hold. Consequently, the corresponding equations of type (\ref{eq0.1}) admit a quadratic first integral. In this connection, it is interesting to find all equations of type (\ref{eq0.1}) which have a quadratic first integral. Therefore, in this section we obtain the necessary and sufficient conditions for the existence of a quadratic first integral for equation (\ref{eq0.1}). We show that these conditions are equivalent to the conditions for an equation of type (\ref{eq0.1}) to possesses a Lax representation with the $L$ matrix of form \eqref{0.3.1}. Therefore, an equation of type (\ref{eq0.1}) admits a Lax representation with the $L$ matrix of form \eqref{0.3.1} if and only if it admits a quadratic first integral and vice versa. First, we consider the generic case, i.e. we do not impose any assumptions on functions $f$ and $g$. Then, we study the case of the classical Li\'enard equation, which is interesting from an applied point of view.

\subsection{The generic case of equation \eqref{eq0.1}}

Suppose that the expression
\begin{equation}
I=(y_{z}+A(z,y))^{2}+B(z,y),
\label{eq:It1}
\end{equation}
where $A$ and $B$ are certain functions, is a first integral for (\ref{eq0.1}). Then
\begin{equation}
\begin{gathered}
DI|_{(\ref{eq0.1})}=2(A_{y}-f)y_{z}^{2}+\cr+(2AA_{y}-2fA+2A_{z}+B_{y}-2g)y_{z}+2AA_{z}-2Ag+B_{z}=0.
\end{gathered}
\label{eq:It3}
\end{equation}
Since functions $A$ and $B$ do not depend on $y_{z}$, after some simplifications, we get  
\begin{equation}
f=A_{y}, \quad g=A_{z}+\frac{B_{y}}{2}, \quad B_{z}-A B_{y}=0.
\label{eq:It5}
\end{equation}
Thus, an equation from family \eqref{eq0.1} admits first integral \eqref{eq:It1} if it is of the form
\begin{equation}
y_{zz}+A_{y} y_{z}+A_{z}+\frac{B_{y}}{2}=0,
\label{eq:It1a}
\end{equation}
provided that $B_{z}-A B_{y}=0$.

Conversely, suppose that the functions $A$ and $B$ satisfy \eqref{eq:It5}. Then, \eqref{eq:It1} is a first integral of \eqref{eq:It1a} since
\begin{equation}
D I=2(y_{z}+A)\left(y_{zz}+A_{y} y_{z}+A_{z}+\frac{B_{y}}{2}\right).
\end{equation}

Therefore, expression \eqref{eq:It1} is a first integral for \eqref{eq0.1} if and only if equation \eqref{eq0.1} is of form \eqref{eq:It1a} and functions $A$ and $B$ satisfy \eqref{eq:It5}. Consequently, compatibility conditions for (\ref{eq:It5}) as an overdetermined system for the functions $A$ and $B$ give us the necessary and sufficient conditions for equation (\ref{eq0.1}) to possess a first integral of form (\ref{eq:It1}).  Moreover, one can show that if we denote $U^{2}$ by $B$ and $F$ by $A$ in (\ref{eq3.1.1.2.1}) we obtain exactly (\ref{eq:It5}). Thus, the existence of first integral (\ref{eq:It1}) of (\ref{eq0.1}) is completely equivalent to the existence of a Lax representation with the $L$ matrix of form (\ref{0.3.1}).

Let us explicitly find correlations on functions $f$ and $g$ that provide compatibility of (\ref{eq:It5}). We introduce the following notations
\begin{equation}
\begin{gathered}
S=g_{yyy}f_{y}+f_{yy}f_{zy}-f_{yy}g_{yy}-f_{zyy}f_{y}, \quad P=f_{z}-g_{y}, \quad Q=g_{yy},  \cr
R=f_{yy}, \quad T=f_{zz}-g_{yz}-2ff_{z}+2f g_{y}.
\end{gathered}
\label{eq:It11}
\end{equation}
Depending on values of $P$, $S$, $T$ and $f_{y}$, the compatibility conditions for (\ref{eq:It5}) split into several cases.

First, we proceed with the most general case. Assume that $S\neq0$ and $f_{y}\neq0$, then we compare various mixed partial derivatives of $B$ with respect to $z$ and $y$ up the fourth order. As a result, we obtain
\begin{equation}
\begin{gathered}
\left( fP_{y}f_{y}-Pf_{y}^{2}+QP_{y}+2 P_{y}^{2}-f_{y}P_{zy} \right)  \left( 2 PRf-3 Pf_{y}^{2}-3 fP_{y}f_{y}- \right. \cr \left. -RP_{z}+f_{y}P_{zy} \right)
-\left( 2 P f^{2}f_{y}-2 P^{2}f_{y}+2 PQf+4 PfP_{y}-2 Pf_{y}g_{y}-\right. \cr \left.-3 fP_{z}f_{y}-gP_{y}f_{y}-QP_{z}-2 P_{z}P_{y}+f_{y}P_{zz} \right) S=0, \vspace{0.1cm}\cr
S^{2}P_{y}+ \left( 4 ff_{y}+Q \right) {S}^{2}-SS_{z}f_{y}+f_{y} \left( 2 Pf-P_{z} \right) R_{y}S-R P_{y}^{2}S+\cr+ \left( f_{y}R_{z}-4 Rff_{y}-6 f_{y}^{3}-QR \right) P_{y}S-f_{y}R \left( 4 Pf_{y}-P_{zy} \right) S+\cr
+3 S_{y}f_{y}^{2}fP_{y}-S_{y}f_{y} \left( 2 PRf-3 Pf_{y}^{2}-RP_{z}+P_{zy}f_{y} \right)=0.
\end{gathered}
\label{eq:It_cc_generic_case}
\end{equation}
The parameters of first integral (\ref{eq:It1}) in this case can be found via the relations
\begin{equation}
\begin{gathered}
A=\frac { \left( P_{z}-2Pf \right) R+3 Pf_{y}^{2}+ \left( 3 fP_{y}-P_{zy} \right) f_{y}}{S}, \cr
B_{y}=\frac {2(2 Pf-P_{{z}})}{f_{y}}-\frac {2P_{y} \left( 2 PRf-3Pf_{y}^{2}-3 fP_{y}f_{y}-RP_{z}+P_{zy}f_{y} \right)}{f_{y}S},\cr
B_{z}=A B_{y}.
\end{gathered}
\label{eq:It_cc_generic_case_FI}
\end{equation}
Here it is assumed that $S\neq0$ and $f_{y}\neq0$.

Let us briefly describe the calculations involved in obtaining (\ref{eq:It_cc_generic_case}) and (\ref{eq:It_cc_generic_case_FI}). Both expressions for computing the functions $A$ and $B$ and conditions on the functions $f$ and $g$ are obtained by considering (\ref{eq:It5}) as an overdetermined system for $A$ and $B$ and finding its compatibility conditions. According to the Riquier--Janet theory (see \cite{Reid1996} and references therein), this is done via cross differentiation of equations from (\ref{eq:It5}) and comparing various mixed partial derivatives of the functions $A$ and $B$. Notice that system (\ref{eq:It5}) is already in the simplified orthonomic form \cite{Reid1996}. The comparison of second and third order mixed derivatives of $B$ does not result in any conditions on functions $f$ and $g$, but gives us expressions for $A_{zz}$ and $A_{z}$. If one proceeds further, taking into account expressions obtained at the previous step, the comparison of $B_{yyzy}$ with $B_{zyyy}$ results in the expression for $A$ given in (\ref{eq:It_cc_generic_case_FI}). Then, we compare $B_{zyyz}$ with $B_{yyzz}$ and $B_{yzyz}$ with $B_{yyzz}$ correspondingly. As a result, we obtain compatibility conditions (\ref{eq:It_cc_generic_case}). Computing further mixed derivatives of $B$ with respect to $z$ and $y$ does not yield any new conditions on the functions $f$ and $g$. The values of $B_{z}$ and $B_{y}$ can be easily obtained with the help of equations (\ref{eq:It5}) and the value of $A$ presented in (\ref{eq:It_cc_generic_case_FI}).

Let us note that there is an alternative way of deriving compatibility conditions (\ref{eq:It_cc_generic_case}). Indeed, by cross differentiating the second and the third equations from (\ref{eq:It5}) one can exclude the function $B$ from (\ref{eq:It5}). As a consequence, we obtain a linear overdetermined system of two equations for the function $A$. Compatibility conditions for this system of equations can be obtained as follows. First, comparison of $A_{yzz}$ with $A_{zzy}$ does not lead to any compatibility conditions, but allows us to find the expression for $A_{z}$. As a result, we can compare $A_{zy}$ with $A_{yz}$ and obtain the expression for $A$ given in (\ref{eq:It_cc_generic_case_FI}). Then, substituting expressions for $A_{z}$ and $A$ into the initial system of equations we obtain the first compatibility condition from (\ref{eq:It_cc_generic_case}). Finally, comparing $A_{zyy}$ and $A_{yyz}$ we obtain the second compatibility condition from (\ref{eq:It_cc_generic_case}). Calculation of further mixed partial derivatives of $A$ does not lead to any new compatibility conditions. Let us remark that we believe that the both ways for obtaining compatibility conditions (\ref{eq:It_cc_generic_case}) and relations (\ref{eq:It_cc_generic_case_FI}) require similar amount of calculations.

In order to verify that we obtain all integrability conditions for (\ref{eq:It5}) we compare our results with the results obtained with the help of the Rif package \cite{Reid1996}, which is based on the Riquier--Janet theory and differential Gr\"{o}bner basis algorithm. Our results coincide with the results produced by the Rif. Thus, we believe that we find all compatibility conditions for (\ref{eq:It5}). Notice that further we do not present details of calculations of compatibility conditions, since they are almost the same as those presented above.


Now we deal with the case of $f_{y}=0$. We assume that $g_{yy}\neq0$, since otherwise we obtain a trivial linear equation. Then, the compatibility conditions are the following
\begin{equation}
\begin{gathered}
2 QQ_{z}+2 TQ_{y}-6 Q^{2}f=0, \quad \left( fg-g_{z} \right) Q^{5}+3 Q_{y}T \left( 3 Tf-T_{z} \right) Q^{2}+\cr+\left(4 Tf_{z}-12 T f^{2}-Tg_{y}+7 fT_{z}-T_{zz} \right) Q^{4}+Q_{yy}T^{3}Q-3 Q_{y}^{2} T^{3}=0.
\end{gathered}
\label{eq:It_cc_f_y=0}
\end{equation}
The parameters of first integral (\ref{eq:It1}) can be found from the relations
\begin{equation}
\begin{gathered}
A=-\frac{T}{g_{yy}},  \quad
B_{y}=\frac {2(gg_{yy}^{2}-2\,Tfg_{yy}+TT_{y}+g_{yy}T_{z})}{g_{yy}^{2}}, \quad B_{z}=A B_{y}.
\end{gathered}
\label{eq:It_cc_f_y=0_FI}
\end{equation}

Further, we proceed with the case of $S=0$. Consequently, we find the following correlations on $f$ and $g$
\begin{equation}
\begin{gathered}
RP_{y}-f_{y}P_{yy}=0, \quad  Pf_{y}^{2}+fP_{y}f_{y}+RT-f_{y}T_{y}=0,\vspace{0.1cm}\cr
 \left( Q_{y}P_{y}-RfP_{y}-TR_{y} \right) f_{y}-6 f_{y}^{3}P_{y}-Rf_{y}^{2}P-R \left( QP_{y}-RT \right)=0,
\end{gathered}
 \label{eq:It_cc_S=0}
\end{equation}
with two additional constraints that are given by
\begin{equation}
\begin{gathered}
Rf_{y}QP_{y}T_{z}-Q_{y}f_{y}QTP_{y}-RQ_{z}f_{y}TP_{y}-f_{y}^{2}T \left( 3{f}^{2}+2P+g_{y} \right) RP_{y}+
\cr
+Rf_{y}^{2}fP_{y}T_{z}+2Q_{y}f_{y}^{2}fTP_{y}-4Rf_{y}^{2}QTP-Q_{y}f_{y}^{2}P_{y}T_{z}-{R}^{2}f_{y}f{T}^{2}+
\cr
+R f_{y}^{3}PT_{z}+R{Q}^{2}TP_{y}+RQT{P_{y}}^{2}-Q_{y}f_{y}T{P_{y}}^{2}-
\cr
-f_{y}^{3}T \left( 5Pf+7T \right) R+3RQ_{y}f_{y}{T}^{2}+2Q_{y}f_{y}^{3}TP-12f_{y}^{3}QTP_{y}-
\cr
-18f_{y}^{4}fTP_{y}+f_{y}^{2}TQ_{zy}P_{y}+6f_{y}^{4}P_{y}T_{z}-2{R}^{2}Q{T}^{2}
\cr
-6f_{y}^{3}TP_{y}^{2}-f_{y}^{2}Q_{yy}T^{2}-18f_{y}^{5}TP-3Rf_{y}QfTP_{y}=0,
\end{gathered}
 \label{eq:It_cc_S=0_ac}
\end{equation}

\begin{equation}
\begin{gathered}
28 f_{y}^{5}gP^{2}-2 P_{y}^{2} \left( 9 PT-9 Tf^{2}+3 Tg_{y}+g_{z}P_{y}\right) f_{y}^{2}-4 T_{zz}f_{y}^{4}P +\cr
+2 P_{y} \left(2 f^{2}gP_{y}+18 PTf-6 PgP_{y}+fP_{y}g_{z}-gP_{y}g_{y}-4 T^{2}-\right. \cr
\left.-3 Tf^{3} \right) f_{y}^{3}- \left( 15 PTf+12 PgP_{y}+4 T^{2} \right) f_{y}^{3}Q+R f_{y}^{2}T_{z}^{2}+\cr
+ \left( 18 Pf+4 T \right) f_{y}^{4}T_{z}+RQ^{2}T^{2}+2 T_{zz}f_{y}^{2}P_{y}^{2}-5 Rf_{y}TP_{y}T_{z}+\cr
+15 Rf_{y}^{3}gTP-2 f_{y}^{2}QfP_{y}T_{z}-2 Rf_{y}QTT_{z}-2 T_{zz}f_{y}^{3}fP_{y}+2 R^{2}f_{y}gT^{2}+ \cr
+TP_{y} \left( 15 Tf-gP_{y} \right) f_{y}R-12 f_{y}^{2}fP_{y}^{2}T_{z}+2 P_{y} \left( 4 f^{2}-P+g_{y} \right) f_{y}^{3}T_{z}+\cr
+4 Q_{z}f_{y}^{3}TP+5 f_{y}^{3}QPT_{z}-5 f_{y}^{2}Q^{2}TP-2 Q_{y}f_{y}T^{2}P_{y}+2RQT^{2}P_{y}-\cr
-P_{y} \left( 3 Tf^{2}+4 fgP_{y}+19 PT+3 Tg_{y}+g_{z}P_{y} \right) f_{y}^{2}Q-\cr
-T \left( Tf+gP_{y} \right) f_{y}^{2}Q_{y}-T_{zz}Rf_{y}^{2}T+Q_{y}f_{y}^{2}TT_{z}-Q_{y}f_{y}QT^{2}+\cr
+RQ_{z}f_{y}T^{2}+P_{y} \left( 3 Tf+gP_{y} \right) Q_{z}f_{y}^{2}-Rf_{y}QgTP_{y}+\cr
+T \left(7 fgP_{y}-2 Tf^{2}+23 PT+Tg_{y}+g_{z}P_{y} \right) f_{y}^{2}R+\cr
+\left( 2 Tf-gP_{y} \right) f_{y}^{2}RT_{z}+T_{zz}f_{y}^{2}QP_{y}-Q_{z}f_{y}^{2}P_{y}T_{z}+3 {R}^{2}T^{3}+\cr
+ \left(20 PfgP_{y} -14 PTf^{2}+44 P^{2}T+4 PTg_{y}+4 PP_{y}g_{z}\right. \cr \left.-4 T^{2}f-4 TgP_{y} \right) f_{y}^{4}=0.
\end{gathered}
 \label{eq:It_cc_S=0_ac_1}
\end{equation}

Let us note that the first condition form (\ref{eq:It_cc_S=0}) is the condition $S=0$, which is presented in terms of $R$ and $P$.

The functions $A$ and $B$ are defined by the relations
\begin{equation}
\begin{gathered}
A=\frac {T_{z}f_{y}-Tff_{y}-gP_{y}f_{y}-QT-2 TP_{y}}{4 Pf_{y}^{2}+2 fP_{y}f_{y}-QP_{y}+RT-2 P_{y}^{2}},\cr
B_{y}= \frac {2(f_{y}P_{y}T_{z}-4 PTf_{y}^{2}-3 TfP_{y}f_{y}-gP_{y}^{2}f_{y}-RT^{2})}{f_{y} \left( 4 P{f_{y}}^{2}+2
 fP_{y}f_{y}-QP_{y}+RT-2 {P_{y}}^{2} \right) }
,\quad
B_{z}= A B_{y}.
\end{gathered}
\label{eq:It_cc_S=0_FI}
\end{equation}
Here we suppose that $4 Pf_{y}^{2}+2 fP_{y}f_{y}-QP_{y}+RT-2 P_{y}^{2}\neq0$, $P_{y}\neq0$, $T\neq0$ and $f_{y}\neq0$.

Let us assume that $S=0$ and $T=0$. As a result, we get that
\begin{equation}
f_{z}-g_{y}=0, \quad A_{z}=g, \quad A_{y}=f, \quad B=\mbox{const},
 \label{eq:It_cc_S=T=0}
\end{equation}
where, without loss of generality, $B$ can be set equal to zero. In order to exclude the case of a trivial linear equation, we assume that $|f_{y}|^2+|g_{yy}|^2\neq0$ in \eqref{eq:It_cc_S=T=0}.

Another separate case, which we need to consider is the case of $T\neq0$ and $R T+4f_{y}^{2}P-2P_{y}^{2}+(2ff_{y}-Q)P_{y}=0$.  This leads to
\begin{equation}
\begin{gathered}
 R T+4f_{y}^{2}P-2P_{y}^{2}+(2ff_{y}-Q)P_{y}=0, \quad RP_{y}-f_{y}P_{yy}=0,\cr
f_{y} \left(QP_{y}-4 P f_{y}^{2}-2 fP_{y}f_{y}+2 P_{y}^{2} \right) R_{y}-R \left(2 RP_{y}^{2}-5 PR f_{y}^{2}-\right. \cr \left.-3 RfP_{y}f_{y} -6 P_{y}f_{y}^{3}+Q_{y}P_{y}f_{y}\right)=0, \cr
gf_{y}{R}^{2}P_{y}+2 Rf_{y}^{3}fP+ \left( 5 PQ+2 PP_{y}+2 g_{y}P_{y} \right) R f_{y}^{2} -8 f_{y}^{4}fP_{y}+\cr
+P_{y} \left( 2 Qf-2 fP_{y}-Q_{z} \right) Rf_{y}-P_{y}^{2} \left( Q+2P_{y} \right)R-16 f_{y}^{5}P+\cr+ \left(4QP_{y}-4 PQ_{y}+8 P_{y}^{2} \right) f_{y}^{3}-2 Q_{y}f_{y}^{2}fP_{y}+Q_{y}P_{y} \left( Q+2 P_{y} \right) f_{y}=0.
\end{gathered}
 \label{eq:It_cc_S=0_1}
\end{equation}
At the same time we find that
\begin{equation}
\begin{gathered}
A_{y}=f, \quad A_{z}=\frac { \left( Rg-2 fP_{y} \right) f_{y}+2 P_{y}^{2}+ \left( Q -AR\right) P_{y}-4 P f_{y}^{2}}{Rf_{y}}
,\cr
B_{y}=\frac {2(ARP_{y}+4 Pf_{y}^{2}+2 fP_{y}f_{y}-QP_{y}-2 P_{y}^{2})}{Rf_{y}}, \quad B_{z}=A B_{y}.
\end{gathered}
 \label{eq:It_cc_S=0_1_FI}
\end{equation}
Here we assume that $R=f_{yy}\neq0$.

Finally, we suppose that $T\neq0$, $R T+4f_{y}^{2}P-2P_{y}^{2}+(2ff_{y}-Q)P_{y}\neq0$ and $P_{y}=f_{yz}-g_{yy}=0$. Consequently, we get
\begin{equation}
\begin{gathered}
f_{yz}-g_{yy}=0, \quad (2f g_{y}-2ff_{z}+f_{zz}-g_{zy})R+3(f_{z}-g_{y})f_{y}^{2}=0,\cr
g R^{3}+(3f_{y}g_{y}+3fQ+3Q_{z})R^{2}-\cr-(27Q f_{y}^{2}+3f f_{y}Q_{y}+3f_{y}Q_{zy}-6QQ_{y})R+6f_{y}Q_{y}(3f_{y}^{2}+Q_{y})=0,\cr
4R^{2}-3f_{y}R_{y}=0,\quad 3f_{y}^{2}Q_{yy}+4Q R^{2}-8R f_{y}Q_{y}=0,
\end{gathered}
 \label{eq:It_cc_S=0_Py=0}
\end{equation}
and
\begin{equation}
A=\frac{3(3f_{y}^{3}-ff_{y}R+Q_{y}f_{y}-QR)}{R^{2}}, \quad B_{y}=2(g-A_{z}), \quad B_{z}=A B_{y},
 \label{eq:It_cc_S=0_Py=0_FI}
\end{equation}
where $R=f_{yy}\neq0$.

Therefore, we have found the necessary and sufficient conditions for equation (\ref{eq0.1}) to possess a first integral in the form (\ref{eq:It1}). These conditions splits into six separate cases: generic case (\ref{eq:It_cc_generic_case}) and special cases (\ref{eq:It_cc_f_y=0}), (\ref{eq:It_cc_S=0}) and (\ref{eq:It_cc_S=0_ac}), (\ref{eq:It_cc_S=0_ac_1}), (\ref{eq:It_cc_S=T=0}), (\ref{eq:It_cc_S=0_1}) and (\ref{eq:It_cc_S=0_Py=0}). Notice that the case when equation \eqref{eq0.1} degenerates into a linear one (that is, when $|f_{y}|^2+|g_{yy}|^2=0$) is not considered in this work, since in the compatibility conditions obtained above it is assumed that $f_{y}\neq0$, $g_{yy}\neq0$ or $|f_{y}|^2+|g_{yy}|^2\neq0$. Below, we consider a special case of equation (\ref{eq0.1}), namely, the Li\'enard equation, which is interesting from an applied point of view.

\subsection{The classical Li\'enard equation}

Let us consider equation (\ref{eq0.2}) and obtain the necessary and sufficient conditions for the existence of first integral (\ref{eq:It1}). Overdetermined system of equations for the parameters of this first integral is completely the same as in the case of equation (\ref{eq0.1}). Thus, we need to obtain the compatibility conditions for (\ref{eq:It5}) taking into account that $f_{z}=g_{z}=0$. This can be easily done as in the case of equation (\ref{eq0.1}). As a result, we get
\begin{equation}
f g g_{yy}-g_{y}(2 f g_{y}+f_{y} g)=0, \quad A=-\frac{fg}{g_{y}}, \quad B_{y}=2g, \quad B_{z}=A B_{y},
\label{eq:Lienard_cc}
\end{equation}
where it is assumed that $g_{y}\neq0$. The case of $g_{y}=0$ leads to a trivial Li\'enard equation. Notice also that $A$ does not depend on $z$. It means that the classical Li\'enard equation can admit a first integral in the form (\ref{eq:It1}) only with $A_{z}=0$.

The condition on functions $f$ and $g$ from (\ref{eq:Lienard_cc}) can be explicitly solved with respect to either function $g$, which yields
\begin{equation}
g=\left(\kappa \int f {\rm d}y +\mu\right)^{-1},
\label{eq:Lienard_cc_explicit}
\end{equation}
or function $f$,
\begin{equation}
f=-\nu g^{-2} g_{y}.
\label{eq:Lienard_cc_explicit_1}
\end{equation}
Here $\kappa\neq0$, $\mu$ and $\nu\neq0$ are arbitrary constants. Either condition (\ref{eq:Lienard_cc_explicit}) or condition (\ref{eq:Lienard_cc_explicit_1}) explicitly defines a family of classical Li\'enard equations, which possess first integral (\ref{eq:It1}).

\section{Examples}

In this section we present several examples of equations from family (\ref{eq0.1}), with generalizations of the Van der Pol and Duffing equations among them, which simultaneously possess a quadratic first integral and a Lax representation.

\textbf{Example 1.} Consider the classical Li\'enard equation (\ref{eq0.2}) with
\begin{equation}\label{eq:e1_1}
g= \alpha(y^2 + 1),
\end{equation}
where $\alpha\neq0$ is an arbitrary parameter and suppose that this equation admits a quadratic first integral. Then, from (\ref{eq:Lienard_cc_explicit_1}) we find that $f=\beta y (y^{2}+1)^{-2}$, where $\beta=2\nu/\alpha\neq0$ is an arbitrary constant. Therefore,
the Li\'enard equation
\begin{equation}\label{eq:e1_2}
y_{zz} +\frac{\beta y}{(y^2+1)^{2}}y_z+\alpha (y^2 + 1) = 0,
\end{equation}
admits the following Lax pair
\begin{equation}\label{eq:e1_3}
  L = \left(\begin{array}{cc}
              y_z - \frac{\beta}{2 (y^2+1)} & U \\
              U & -y_z +\frac{\beta}{2 (y^2+1)}
            \end{array}\right),\quad M=\left(\begin{array}{cc}
                                                0 & \frac{U_y}{2}\\
                                                  -\frac{U_y}{2}&  0
                                             \end{array}\right),
\end{equation}
where
\begin{equation}
  U^2=\frac{2\alpha}{3} y^3+ 2\alpha y -\alpha\beta z,
\end{equation}
and the following first integral
\begin{equation}\label{eq:e1_4}
I = \left(y_z - \frac{\beta}{2(y^2+1)}\right)^2 + \frac{2\alpha}{3}y^3 + 2\alpha y -\alpha \beta z.
\end{equation}
Notice that equation (\ref{eq:e1_2}) can be considered as an anharmonic oscillator with nonlinear friction.

\textbf{Example 2.} Suppose that $f=z/(2y+z)^2$ and $g=\alpha^2(2y^3+3zy^2+z^2y)-y/(2y+z)^2$, where $\alpha\neq0$ is an arbitrary constant. One can show that these functions $f$ and $g$ satisfy conditions (\ref{eq:It_cc_generic_case}) and, thus, the corresponding equation of type (\ref{eq0.1}) admits a Lax representation and a quadratic first integral. Indeed, the equation
\begin{equation}\label{3.4}
  y_{zz} + \frac{z}{(2y+z)^2}y_z + \alpha^2(2y^3+3zy^2+z^2y)-\frac{y}{(2y+z)^2}=0,
\end{equation}
has the Lax pair
\begin{equation}\label{eq3.5}
  L=\left(\begin{array}{cc}
            y_z + \frac{y}{2y+z}  & \alpha(y^2+zy) \\
           \alpha(y^2+zy) & - y_z - \frac{y}{2y+z}
          \end{array}\right), \quad
          M = \left(\begin{array}{cc}
                0 & \frac{\alpha(2y+z)}{2} \\
                - \frac{\alpha(2y+z)}{2}  & 0
              \end{array}\right),
\end{equation}
and the first integral
\begin{equation}\label{3.6}
I = \left(y_z+  \frac{y}{2y+z}\right)^{2} + \alpha^2y^2(y+z)^2.
\end{equation}

\textbf{Example 3.} Now we demonstrate an example of an equation of type (\ref{eq0.1}), for which we cannot use conditions (\ref{eq:It_cc_generic_case}) since in this case $S=0$.  Suppose that $f={\rm e}^{-\alpha z}(y+\delta)^{-2}$ and $g=\alpha {\rm e}^{-\alpha z}(y+\delta)^{-1}+y+\delta$, then, one can show that these values of $f$ and $g$ satisfy conditions (\ref{eq:It_cc_S=0_Py=0}). As a result, we find that the Lax pair
\begin{equation}
\begin{gathered}
  L=\left(\begin{array}{cc}
      y_{z}-\frac{{\rm e}^{-\alpha z}}{y+\delta} & \left(y^{2}+2\delta y+\frac{2}{\alpha}{\rm e}^{-\alpha z}\right)^{\frac{1}{2}} \\
      \left(y^{2}+2\delta y+\frac{2}{\alpha}{\rm e}^{-\alpha z}\right)^{\frac{1}{2}} & -y_{z}+\frac{{\rm e}^{-\alpha z}}{y+\delta}
    \end{array}\right),  \hfill \quad \quad \quad \quad \quad \quad \quad \vspace{0.3cm}\cr
               \quad \quad \quad \quad \quad \quad \quad M=\left(\begin{array}{cc}
                                         0& \frac{y+\delta}{2\left(y^{2}+2\delta y+\frac{2}{\alpha}{\rm e}^{-\alpha z}\right)^{\frac{1}{2}}} \\
                                      -\frac{y+\delta}{2\left(y^{2}+2\delta y+\frac{2}{\alpha}{\rm e}^{-\alpha z}\right)^{\frac{1}{2}}} & 0
                                    \end{array}\right),
\end{gathered}
\end{equation}
and the first integral
\begin{equation}
 I = \left(y_{z}-\frac{{\rm e}^{-\alpha z}}{y+\delta}\right)^{2}+y^{2}+2\delta y+\frac{2}{\alpha}{\rm e}^{-\alpha z},
\end{equation}
correspond to the equation
\begin{equation}
y_{zz}+\frac{{\rm e}^{-\alpha z}}{(y+\delta)^{2}} y_{z}+\frac{\alpha {\rm e}^{-\alpha z}}{(y+\delta)}+y+\delta=0.
\label{eq:ex3_1}
\end{equation}
Let us note that equation \eqref{eq:ex3_1} can be considered as a generalization of a harmonic oscillator. Suppose that $\alpha>0$, then one can see that for the large values of either $z$ or $y$, equation \eqref{eq:ex3_1} tends to a harmonic oscillator.

\textbf{Example 4.} Consider an equation of type (\ref{eq0.1}) with $f_{z}=0$. The functions $f=-\alpha(1-y^{2})$ and $g=\beta y {\rm e}^{-2\alpha z} (y^{2}-3)^{-2}$, where $\alpha\neq0$ and $\beta\neq0$ are arbitrary parameters, satisfy conditions  \eqref{eq:It_cc_generic_case}. Therefore, the equation
\begin{equation}
  y_{zz}-\alpha(1-y^{2})y_{z}+\frac{\beta y}{ {\rm e}^{2\alpha z} (y^{2}-3)^{2}}=0,
  \label{eq:e5_1}
\end{equation}
has the Lax pair
\begin{equation}
\begin{gathered}
L=\left(\begin{array}{cc}
          y_z +\frac{\alpha y}{3}(y^2-3) & \sqrt{\frac{-\beta}{3(y^2-3)}}y{\rm e}^{-\alpha z} \\
          \sqrt{\frac{-\beta}{3(y^2-3)}}y{\rm e}^{-\alpha z} & -y_z -\frac{\alpha y}{3}(y^2-3)
        \end{array}\right),\cr M=\left(\begin{array}{cc}
                                            0 & -\frac{\sqrt{-3\beta}{\rm e}^{-\alpha z}}{2(y^2-3)^{\frac{3}{2}}} \\
                                           \frac{\sqrt{-3\beta}{\rm e}^{-\alpha z}}{2(y^2-3)^{\frac{3}{2}}} &  0
                                         \end{array}\right),
\end{gathered}
\end{equation}
and the first integral
\begin{equation}
 I = \left(y_{z}+\frac{\alpha y(y^{2}-3)}{3}\right)^{2}-\frac{\beta y^2 {\rm e}^{-2 \alpha z}}{3(y^2-3)}.
\end{equation}
Equation (\ref{eq:e5_1}) can be considered as a non-autonomous generalization of the Van der Pol oscillator.

\textbf{Example 5.} Now we demonstrate an example of an equation of type (\ref{eq0.2}). Choosing the Van der Pol like friction $f=-\alpha(1-y^2)$ and using (\ref{eq:Lienard_cc_explicit}) at $\kappa=(\alpha \beta)^{-1}$ and $\mu=0$, where $\alpha\neq0$ and $\beta\neq0$ are arbitrary parameters, we obtain that the equation
\begin{equation}
  y_{zz} - \alpha(1-y^2) y_{z} + \frac{3\beta}{y (y^2 -3)} = 0,
\end{equation}
has the Lax pair
\begin{equation}
  L = \left(\begin{array}{cc}
              y_z - \alpha y(1-\frac{y^2}{3})& U \\
              U & -y_z + \alpha y(1-\frac{y^2}{3})
            \end{array}\right), \quad M=\left(\begin{array}{cc}
                                                0 & U_{y}/2 \\
                                                -U_{y}/2 &  0
                                              \end{array}\right),
\end{equation}
with
\begin{equation}
  U = \left[2\alpha\beta z - \beta \ln\left(\frac{y^2}{y^2-3}\right)\right]^{\frac{1}{2}},
\end{equation}
and the first integral
\begin{equation}
  I = \left[y_z - \alpha y\left(1-\frac{y^2}{3}\right)\right]^{2} + 2\alpha\beta z - \beta \ln\left(\frac{y^2}{y^2-3}\right).
\end{equation}
This example can be considered as a generalization of the Van der Pol equation, which has a quadratic first integral.

\textbf{Example 6.} Finally we demonstrate an example of an equation of type (\ref{eq0.2}), which generalizes the Duffing equation. Let $g=y^3+\alpha y + \beta$, then using (\ref{2.5}) we obtain that the equation
\begin{equation}
  y_{zz} -\frac{\nu(3y^2+\alpha)}{(y^3+\alpha y + \beta)^2}y_z + y^3+\alpha y + \beta = 0,
\label{eq:ex6_1}
\end{equation}
has the Lax pair
\begin{equation}
  L = \left(\begin{array}{cc}
              y_z + \frac{\nu}{(y^3+\alpha y + \beta)} & U \\
              U & -y_z - \frac{\nu}{(y^3+\alpha y + \beta)}
            \end{array}\right),\quad M=\left(\begin{array}{cc}
                                               0 & \frac{g}{2U} \\
                                               -\frac{g}{2U} & 0
                                             \end{array}\right),
\end{equation}
where
\begin{equation}
  U = \left(\frac{y^4}{2}+\alpha y^2 + 2\beta y + 2\nu z\right)^{\frac{1}{2}}.
\end{equation}
Equation \eqref{eq:ex6_1} also admits the first integral
\begin{equation}
  I = \left(y_z + \frac{\nu}{(y^3+\alpha y + \beta)}\right)^2 + \frac{y^4}{2}+\alpha y^2 + 2\beta y + 2\nu z.
\label{eq:ex6_3}
\end{equation}
This example can be considered as the Duffing equation with nonlinear friction. First integral \eqref{eq:ex6_1} can be interpreted as the total energy of the considered system, where the term $\nu z$ is responsible for friction, while the rest of the terms represent kinetic and potential energy. Notice that at $\nu\rightarrow 0$ from \eqref{eq:ex6_1} and \eqref{eq:ex6_3} we obtain the classical Duffing equation and its first integral, correspondingly.



\section{Conclusion and discussion}
Here we summarize the results obtained in this work. First, it is easy to see that if we denote $U^{2}$ by $B$ and $F$ by $A$ in (\ref{eq3.1.1.2.1}) we obtain exactly (\ref{eq:It5}). Thus, sufficient conditions for the existence of a Lax representation with the $L$ matrix of form \eqref{0.3.1} coincide with the necessary and sufficient conditions for the existence of quadratic first integral (\ref{eq:It1}).  Therefore, the following proposition holds:

\textbf{Proposition 1.} The following statements are equivalent
\begin{enumerate}[(I)]
  \item equation (\ref{eq0.1}) has a quadratic first integral in the form (\ref{eq:It1}),
  \item equation (\ref{eq0.1}) admits Lax representation (\ref{eq3.3}),
  \item functions $f$ and $g$ from (\ref{eq0.1}) satisfies one of the following conditions: 1) (\ref{eq:It_cc_generic_case}); 2) (\ref{eq:It_cc_f_y=0}); 3) (\ref{eq:It_cc_S=0}), (\ref{eq:It_cc_S=0_ac}) and (\ref{eq:It_cc_S=0_ac_1}); 4) (\ref{eq:It_cc_S=T=0}); 5) (\ref{eq:It_cc_S=0_1}); 6) (\ref{eq:It_cc_S=0_Py=0}).
\end{enumerate}

From this proposition one can obtain the following corollaries:

\textbf{Corollary 1.} The classical Li\'enard equation (\ref{eq0.2}) can possess a first integral of type (\ref{eq:It1}) only with $A_{z}=0$.

\textbf{Corollary 2.} An equation of type (\ref{eq0.1}) admits an autonomous first integral of type (\ref{eq:It1}) if and only if either $g_{z}=0$ and $f=0$ or $f_{z}=0$ and $g=0$.

While Corollary 1 follows from (\ref{eq:Lienard_cc}), to prove the second one we need to assume that $A_{z}=B_{z}=0$ in (\ref{eq:It5}), which immediately leads to corresponding statement.

Notice that first integral (\ref{eq:It1}) can be considered as an analog of the energy conservation law for the dissipative systems. Indeed, if we expand the expression in brackets in (\ref{eq:It1}), the first term may be interpreted as the kinetic energy, the term with first derivative as a friction term and the rest terms as a potential. Therefore, we have obtained that there is a direct correspondence between the existence of an energy-like first integral for a dissipative system from family (\ref{eq0.1}) and the existence of a Lax representation with the $L$ matrix of form \eqref{0.3.1}. This correspondence may be considered as an analogue for the dissipative systems of the direct connection between the Lax and Arnold--Liouville integrability for Hamiltonian systems (see \cite{Babelon1990}).

It is also interesting to compare the results of this work with recent results (see, \cite{Kudryashov2016,Kudryashov2016a,Kudryashov2017a,Sinelshchikov2018,Ruiz2018}) devoted to integrability of the classical Li\'enard equation \eqref{eq0.2} and its generalizations. Since integrability conditions or conditions for the existence of first integrals obtained in these works are given in terms of relations between the functions $f$ and $g$, one can readily compare them with either \eqref{eq:Lienard_cc_explicit} or \eqref{eq:Lienard_cc_explicit_1}. As a result, one can see that neither \eqref{eq:Lienard_cc_explicit} nor \eqref{eq:Lienard_cc_explicit_1} coincide with any of the integrability conditions obtained in \cite{Kudryashov2016,Kudryashov2016a,Kudryashov2017a,Ruiz2018}. Thus, correlation \eqref{eq:Lienard_cc_explicit} (or \eqref{eq:Lienard_cc_explicit_1}) provides a new condition for the existence of a quadratic first integral for \eqref{eq0.2}. What is more, condition \eqref{eq:Lienard_cc_explicit} is different from the condition for linearizability of \eqref{eq0.2} via the generalized Sundman transformations, which also means that equation \eqref{eq0.2} with coefficients satisfying \eqref{eq:Lienard_cc_explicit} cannot be completely integrated via Lie point symmetries (see \cite{Kudryashov2016}). Finally, in work \cite{Sinelshchikov2018} some integrability conditions for \eqref{eq0.1} were found. However, it can be shown that these conditions are again different from correlations on the functions $f$ and $g$ obtained in Section 3.



Let us briefly summarize the main results of this work. We have found all equations of type (\ref{eq0.1}) that simultaneously admit quadratic first integral (\ref{eq:It1}) and Lax representation (\ref{eq0.4}) with the $L$ matrix of form \eqref{0.3.1}. We have shown that these two notions are equivalent for any equation from family (\ref{eq0.1}). We have considered a particular case of equation (\ref{eq0.1}), namely the classical Li\'enard equation, and have explicitly found the corresponding quadratic first integrals and Lax pairs. We believe that this is the first time when the Lax integrability and the existence of quadratic first integrals have been studied for dissipative differential equations of type (\ref{eq0.1}).

\section*{Acknowledgments}
Authors are grateful to anonymous referees and Associate Editor for their valuable comments and suggestions.

This work is supported by Russian Science Foundation, grant number 18--11--00209.

\end{document}